\overfullrule=0pt   \magnification=\magstep1  \input amssym.def
\font\huge=cmr10 scaled \magstep2
\def\i{{\rm i}}  \def\si{\sigma}   \def\eps{\epsilon} \def\la{\lambda}
\def\Z{{\Bbb Z}}  \def\Q{{\Bbb Q}}  \def\C{{\Bbb C}} \def\i{{\rm i}}
\def\qed{\quad\vrule height 6pt width 6pt depth 0pt}

\font\smal=cmr7
\font\smit=cmmi7
\font\smcap=cmcsc10 
{\nopagenumbers \rightline{{September 13, 1999}}

\rightline{{Preliminary Version}}\bigskip 
\bigskip

\centerline{{\huge Congruence Subgroups and}}\medskip
\centerline{{\huge  Rational Conformal Field
Theory}\footnote{$^\dagger$}{{\smal The published version of this paper will
assume slightly more mathematical sophistication; both versions have
equivalent content but this one is a little more pedagogical.}}}\bigskip
\bigskip 
%%%%%%%%

\centerline{{ Antoine Coste}}
\medskip\centerline{{\it CNRS Laboratory of Theoretical Physics,}} 
\centerline{{\it building 210, Paris XI University }} 
\centerline{{\it 91405 Orsay cedex, France }} \smallskip
\centerline{{ Antoine.Coste@th.u-psud.fr }} \bigskip

\centerline{{Terry Gannon}}
\medskip\centerline{{\it Department of Mathematical Sciences,}}
\centerline{{\it University of Alberta}}
\centerline{{\it Edmonton, Canada, T6G 2G1}}\smallskip
\centerline{{tgannon@math.ualberta.ca}}

\bigskip\bigskip 

\centerline{{\bf Abstract}}\medskip

We address here the question of whether the characters of an RCFT are
modular functions for some level $N$, i.e.\ whether the representation of the
 modular group SL$_2(\Z)$ coming from any RCFT is trivial on some
 congruence subgroup.
We prove that if the matrix $T$, associated to $\left(\matrix{1&1\cr 0&1}
\right)\in{\rm SL}_2(\Z)$, has {\it odd} order, then this must be so. When
the order of $T$ is even, we present a simple test which if satisfied
--- and we conjecture it always will be --- implies that the characters for
that RCFT will also be level $N$. We use this to explain three
curious observations in RCFT made by various authors. \medskip

This is the presubmission copy. We are interested in receiving any feedback.
\vfill\eject}\pageno=1

{{\bf 1. Introduction}}\medskip

Associated to a rational conformal field theory (RCFT),  or related structures
such as affine Kac-Moody algebras or rational vertex 
operator algebras (VOAs), is a finite-dimensional representation $\rho$ of 
the (homogeneous) modular group ${\rm SL}_2(\Z)$. In 
particular, we write\footnote{$^1$}{{\smal Note that our choice of
{\smit S}\ is slightly different from that made by some other
authors. This is discussed more fully four paragraphs into section 2.}} 
 $S:=\rho\left(\matrix{0&1\cr -1&0}\right)$ and
$T:=\rho\left(\matrix{1&1\cr 0&1}\right)$. $S$ determines the fusion 
coefficients $N_{ab}^c$ in the RCFT, by Verlinde's formula:
$$N_{ab}^c=\sum_{d\in\Phi}{S_{ad}\,S_{bd}\,S_{cd}^*\over S_{0d}}\eqno(1)$$
where $a,b,c,d\in\Phi$ label the finitely many primary fields. This modular 
representation is realised by
the characters ch$_a(\tau)$, $a\in\Phi$, of the RCFT:
$$\eqalignno{{\rm ch}_a(-1/\tau)&\,=\sum_{b\in\Phi}S_{ab}\,{\rm ch}_b(\tau)&(2a)\cr
{\rm ch}_a(\tau+1)&\,=\sum_{b\in\Phi}T_{ab}\,{\rm ch}_b(\tau)&(2b)\cr}$$
Basic known properties of $\rho$ (i.e.\ $S$ and $T$) will be quickly reviewed
in \S2.

In this paper we address one of the most fundamental questions about 
 these modular functions ch$_a(\tau)$: are they fixed by a congruence
subgroup? The main results of this paper are Theorems 2 and 4. Other
interesting  results are the group presentations in Lemma 1, the
 explanations in \S3, and the two propositions.

Let $\Gamma$ denote SL$_2(\Z)$.
By a {\it congruence subgroup} we mean any subgroup of $\Gamma$ containing
$$\Gamma(N):=\{M\in\Gamma\,|\,M\equiv\left(\matrix{1&0\cr 0&1}\right)
\ ({\rm mod}\ N)\}$$
for some $N$. $\Gamma(N)$ is called the principal congruence subgroup of
level $N$. We are interested in whether a given representation $\rho$ of
$\Gamma$ `factors through a congruence subgroup', i.e.\ whether $\rho$ sends
$\Gamma(N)$ to the identity matrix $I$. This would tell us that $\rho$
is actually a representation of the finite group  
$${\rm SL}_2(N):={\rm SL}_2(\Z/N\Z)\cong\Gamma/\Gamma(N)\ ,$$ consisting 
of all $2\times 2$ matrices $M$ with entries from the integers mod $N$ and 
with determinant
$|M|\equiv 1$ (mod $N$). In other words, the matrices $S$ and $T$ would generate
a group $\langle S,T\rangle$ isomorphic to some factor group of SL$_2(N)$. We 
will say a given
RCFT {\it has the $\Gamma(N)$ congruence property} if its modular representation
$\rho$ is trivial on some $\Gamma(N)$. This implies that its
characters ch$_a$ are all fixed by, i.e.\ modular for,
$\Gamma(N)$. This would mean that the ch$_a$ are all level $N$ modular functions.

See [AS,J] and references therein for some samples from the theory of 
 noncongruence subgroups. The main reason congruence subgroups are so familiar 
is that they arise wherever theta functions of quadratic forms (or lattices)
 do, as was first shown by Hecke and Schoeneberg ({\it c.}\ 1940) --- see
 for instance Chapter VI of [O] for details (our methods result in a
 new proof, given in section 4 below).
In actual fact, {\it congruence subgroups are
far rarer than noncongruence ones}: if we look at an arbitrary subgroup ${\cal 
G}$ of $\Gamma$
with finite but large index $\|\Gamma/{\cal G}\|$, the probability will be almost 1
that it is noncongruence. More precisely [J], the number of noncongruence
 subgroups ${\cal G}$ of index $n=\|\Gamma/{\cal G}\|$ grows faster
 than $({n\over e})^{n\over 6}$, while the number of congruence
 subgroups of index $n$ is bounded above by the much smaller number
 $n^{1+9\,{\rm log}_2n}$.

Or for another indication of their comparative numbers, recall that
given a subgroup ${\cal G}$ of $\Gamma$ with finite index, we can
construct the Riemann surface (with finitely many punctures)
${\cal G}\backslash{\Bbb H}$ from the upper half-plane ${\Bbb H}$. By the {\it
genus} of ${\cal G}$ we mean the genus of that Riemann surface. Then
for any given genus, there are only finitely many congruence subgroups
${\cal G}$ but infinitely many noncongruence subgroups [J].

SL$_2(\Z)$ is truly exceptional. By comparison, all finite-index
subgroups of SL$_n(\Z)$, for $n\ge 3$, are congruence!

The first example of a noncongruence subgroup goes back to Klein
(1879), and we can obtain infinitely many examples as follows.
Consider the function $\xi(\tau):=\eta(\tau)/\eta(13\tau)$, where
$\eta$ is the Dedekind eta. Then $\xi$ is a genus-zero modular function for $\Gamma(26)$,
but for any $m=2,3,4,\ldots$, its $m$th root $\xi(\tau)^{{1\over m}}$ 
(taking the principal branch of log$\,\xi$) is a genus-zero modular function
for a {\it noncongruence} subgroup [AS].

Nevertheless, several people (e.g.\ [Mo,E,ES,DM,BCIR,B]) have conjectured (or at 
least speculated on the possibility of) the following:

\medskip{\smcap Conjecture 1.} {\it All RCFTs have the congruence property, so
in particular 
their characters ch$_a(\tau)$ are modular functions for some} $\Gamma(N)$.\medskip

We will strengthen Conjecture 1 slightly, in \S3.

Why is this conjecture not simply naive optimism? After all, it was
not even known (see below) whether the subgroup of $\Gamma$ fixing the
ch$_a$ has finite index. A reason for suspecting the truth of the
conjecture  is that all known RCFTs possess the congruence property
--- see  for instance 
\S4. But the best motivation for Conjecture 1 is the following 
hope, originally observed empirically by Atkin and Swinnerton-Dyer [AS]:

\medskip{\smcap Conjecture 2.} {\it  
Let $f(\tau)= q^c \   \sum_{n=0}^\infty a_n\,q^{n/b}\not\equiv 0$ be a modular 
function for some subgroup ${\cal G}$ of $\Gamma$ and some $b\in{\Bbb N}$, $c\in\Q$.
If the Fourier coefficients $a_n$ are all algebraic integers, then
${\cal G}$ is a congruence subgroup.}\medskip

The most important examples of algebraic integers here are the
`rational integers' $\Z$, and  the `cyclotomic integers' $\Z[\xi_n]$ given by
polynomials with coefficients in $\Z$, evaluated at some $n$th root of unity
$\xi_n$.

The converse of Conjecture 2 is known to be true (see for example Ch.\ 6 of [L]): the
modular functions for $\Gamma(N)$ with Fourier coefficients
$a_k\in\Z[\xi_N]$ span the space of all modular functions for $\Gamma(N)$.

Now, the coefficients of our RCFT characters
ch$_a$ are in fact rational {integers}, so Conjecture 2 would imply Conjecture 1, at least
if the RCFT characters were linearly independent (which in general they
aren't --- we will return to this important point in \S2).
Thus Conjecture 2 strongly suggests (but in general won't imply)
Conjecture 1.
Although Conjecture 2 seems plausible and would be an important result
in automorphic function theory, it remains unproven.

Why should we care about Conjecture 1? For one thing, congruence subgroups are
much  more familiar, and are also
much better understood --- relatively little is known generally about
{\it non}congruence subgroups. The main technical difficulty with the
latter is the lack for them of a satisfactory theory of Hecke
operators. Whenever the rich theory of modular functions is applied to
RCFT theory, a simplifying (and for many purposes necessary) assumption
is certainly that the characters be modular for a congruence subgroup.
See e.g.\ [ES] for such results, and [E] anticipates that
the classification of all (not necessarily
unitary) RCFTs with effective central charge $\tilde{c}\le 1$, would follow
quickly from some technical modular function results (Serre-Stark)
and the congruence property. It would also be very useful to know that $\rho$
is  in fact a representation of the finite group
SL$_2(N)$, as those groups are so well-understood.
For instance, their representations have been
classified [TNW], and [E]  suggested we should use that to classify the possible
modular data in RCFTs obeying the congruence property.

By comparison, any finite group generated by an order 2 and an order 3
element together --- e.g.\ the alternating group ${\frak A}_n$ for $n\ge 9$ ---
will be a factor group of $\Gamma$ (this is because PSL$_2(\Z)$ is
isomorphic to the free product $\Z/2\Z* \Z/3\Z$). So we can't
expect many interesting general
results on the finite quotients $\Gamma/{\cal G}$, unless we assume in
addition that e.g.\ ${\cal G}$ is a congruence subgroup. 

In that sense, the main value of our paper could be to unlock a door
behind which could lie some mathematical riches for RCFT.

Also, Conjecture 1 would help explain some curiousities in RCFTs (see
 \S3 below).  An example is the observation
made in [BI] that the commutant for $A_\ell^{(1)}$ level $k$ has an integral basis.
More important is the Galois action for any RCFT [CG]: in \S3 we will 
 interpret this 
as the natural Galois action on level $N$ modular functions [L]. 

To help put Conjecture 1 into perspective, consider the much weaker
statement that the matrices $S,T$ for any RCFT must necessarily generate
a finite group (i.e.\ that the subgroup of $\Gamma$ fixing all the
characters ch$_a(\tau)$ has finite index in $\Gamma$). Even that 
was not  known
to be true (but see Theorem 4 below). It is tempting to suspect that
the obvious relations between $S,T$
 stated below in the second paragraph of \S2 
are enough to guarantee this. However, that is a false hope: in actual fact,
the group defined by the presentation\footnote{$^2$}{{\smal To the
left of the bar are the generators, to the right are the relations.}}
$$\langle S,T\,|\,T^N=S^4=(ST)^3=I, S^2\ {\rm and}\ T\ {\rm commute} \rangle$$
equals SL$_2(N)$ for $N\le 5$ but is infinite for all $N>5$. Imposing the
additional condition $S^2=I$ doesn't help: $\langle S,T\rangle$ then equals 
[CM] the symmetric group ${\frak S}_3$
for $N=2$, the tetrahedral group ${\frak A}_4$ for $N=3$, the octahedral
group ${\frak S}_4$ for $N=4$, the icosahedral group ${\frak A}_5$
for $N=5$, and again is infinite for any $N>5$. 
This $\langle S,T\rangle$ (with $S^2=I$) is called the triangle (or
polyhedral) group $(2,3,N)$. 

Thus for $\langle S,T\rangle$ to be finite, we need additional `nonobvious'
relations between the matrices $S$ and $T$. (Imposing additional
relations is equivalent to quotienting $(2,3,N)$ by some normal
subgroup.) At least for
$N$ not a multiple of 6, we find in Theorem 2 exactly one additional
relation that will accomplish this; moreover, it is easy to check this
relation in practice and surprisingly when it holds (which we believe is
always) it tells us this finite group will be a factor group of
SL$_2(N)$.  Indeed, in Theorem
4 we apply this test to prove Conjecture 1 for $N$ odd. 

Incidentally, all finite factor groups of the triangle group $(2,3,6$)
have been classified by Newman (1964). The finite factor groups of
(2,3,7) are important in Riemann surface theory 
and are called Hurwitz groups.

In the next section we give a simple test for determining whether or not a 
given
RCFT has the congruence property, and we prove the congruence property when
$N$ is odd. While our test falls just short of establishing
that all RCFTs (i.e.\ also $N$ even) must have the congruence property, it 
demonstrates why generic
RCFTs should. In the process we obtain natural presentations of the
group SL$_2(N)$. This simplifies 
the congruence subgroup test of [H], as we discuss briefly in \S4.

\bigskip
{{\bf 2. When an RCFT has the congruence property}}\medskip

{\it 2.1. RCFT modular data.}
Consider the matrices $S,T$ corresponding to a given RCFT. We will explicitly
state in the next two paragraphs all properties of $S,T$ we will need.

Because $S$ and $T$ correspond to a representation of $\Gamma$, we know
$(ST)^3=I$, and the {\it charge-conjugation} matrix $C:=S^2$ commutes with both
$S$ and $T$. The matrices $S,T$ are both unitary and symmetric, and
$T$ is  diagonal.
Also, we know [AM] there exists some integer $N>0$ for which
 $T^N=I$. Automatically in RCFT the characters are holomorphic in the upper
half-plane.  Incidentally, this also implies from (2) that each ch$_a$ will
be `meromorphic at each cusp' --- e.g.\ meromorphicity at the cusp
$\i\infty$ would mean that each ch$_a$ has a Laurent expansion in the
local coordinate $q^{{1\over N}}$.
Hence the $\Gamma(N)$ congruence property would imply the characters
are all level $N$ modular functions (the converse is not necessarily
true, because the characters will be linearly dependent in general ---
see below).

We know from [CG] that the entries $S_{ab}$ of $S$ must lie
in a cyclotomic extension ${\Bbb Q}[\xi_n]$ of ${\Bbb Q}$. $\xi_n$ here
is the $n$th root of unity $\exp[2\pi\i/n]$, and $\Q[\xi_n]$ can be thought
of as all complex numbers of the form $a_0+a_1\xi_n+\cdots+a_k\xi_n^k$,
where the coefficients $a_i$ are rational. The Galois group Gal$(\Q[\xi_n]/\Q)$
of ${\Bbb Q}[\xi_n]$ is defined to be the automorphisms of the field
$\Q[\xi_n]$ which fix $\Q$. This Galois group is isomorphic to the 
multiplicative (mod $n$) group of integers
coprime to $n$, which we write $\Z_n^*  $: specifically, the Galois 
automorphism $\si_\ell$
corresponding to $\ell\in\Z_n^* $ takes the number
$a_0+a_1\xi_n+\cdots+a_k\xi_n^k$ to $a_0+a_1\xi_n^\ell+
\cdots+a_k\xi_n^{\ell k}$. 
Think of $\si_\ell$ as a generalisation of complex conjugation --- in fact
complex conjugation equals $\si_{-1}$. Now,
choose any $\si\in\Z_n^* $, then [CG]
$$\si(S_{ab})=\eps_\si(a)\,S_{\si a,b}=\eps_\si (b)\,S_{a,\si b}\eqno(3a)$$   
where $\eps_\si(a)\in\{\pm 1\}$ are signs, and $a\mapsto \si a$ defines a permutation
of $\Phi$, independent of $b$. 
 If we define for each such $\si$ the matrix $G_\si$ given by
$(G_\si)_{ab}=\eps_\si(a)\,\delta_{b,\si(a)}$, then the assignment $\si\mapsto
G_\si$ defines a representation of $\Z_n^* $ and (3a) reads
$$\si(S)=G_\si S=SG_\si^{-1}\ .\eqno(3b)$$
For instance, charge-conjugation $C=G_{-1}$.
This important Galois action (3a) holds for any RCFT, and is a consequence of the 
basic properties of $S,T$ given in the previous paragraph, together with the
fact that the fusion coefficients $N_{ab}^c$ in (1) are rational.
All of these are basic ingredients in any
 RCFT [MS].

A minor clarification should be made. There is an equally valid alternate
choice for $S$, namely $\rho\left(\matrix{0&-1\cr 1&0}\right)$, which is
more commonly made in the literature, resulting in slightly different formulas.
These two possibilities for $S$ are complex conjugates of each other.
For example our
equations $(ST)^3=(TS)^3=I$ would become $(ST)^3=(TS)^3=C$. The choice we
have made seems to result in slightly cleaner formulas.

Strictly speaking, the RCFT characters ch$_a(\tau)$ will not in general be
linearly independent and so equation (2a) will not uniquely determine $S$.
The simplest example of this (there are others) is that $a$ and its
charge-conjugate $Ca$ will always have equal characters ch$_a={\rm ch}_{Ca}$,
even though often $a\ne Ca$. 
The obvious way out is to introduce additional variables in addition to $\tau$,
so ch$_a$ then would involve a more `sensitive' trace. The most familiar 
instance is 
the transition from the theta function $\theta(\tau)=\sum_{n\in\Z}q^{n^2/2}$ 
to
its Jacobi form $\theta(\tau,z)=\sum_{n\in\Z}q^{n^2/2}r^n$, where $q=e^{2\pi \i
\tau}$ and $r=e^{2\pi\i z}$. This is precisely what is done for the affine 
Kac-Moody algebras [KP]. This technical point unfortunately is usually 
overlooked
in the literature, and we will return to it in \S3 (see also [GG]). Until it gets clarified
though, and we learn how to obtain $S$ unambiguously from the RCFT characters,
it will be much more difficult to apply modular (or Jacobi) function theory
rigourously and nontrivially to RCFT.

\smallskip{\it 2.2. A natural presentation of $SL_2(N)$.}
A key tool we need are  generators
and relations for the finite group SL$_2(N)$. There has been quite an
industry in this direction (see e.g.\ [M,BM,CR,H,H2]).
 With a little effort we can write them in
 the following more transparent  form.

\medskip{\smcap Lemma 1.} {\it Choose any $N\in\{1,2,3,\ldots\}$.
By `${1\over p}$' we mean the integer-valued multiplicative inverse of $p$ 
mod $N$. Then we get the following presentations for SL$_2(N)$:

\item{(a)} For $N$ coprime to $p$, where $p$ either equals 2 or 3,
$${\rm SL}_2(N)=\langle s,t\,|\,t^N=s^4=1,(st^{-1})^3=s^2,\,
gs=sg^{-1},\, gt=t^{p^2}g\ {\rm where}\ g:=st^{{1\over p}}st^pst^{{1\over
p}} \rangle$$

\item{(b)} For $N$ coprime to $p$, where $p$ either equals 5 or 7, 
$$\eqalign{{\rm SL}_2(N)=\langle s,t\,|\,&t^N=s^4=1,(st^{-1})^3=s^2,
gs=sg^{-1}, gt=t^{p^2}g, g=t^{p{p-1\over 2}}st^{-{2\over
p}}st^{-{p-1\over 2}}st^2s,\cr & {\rm where}\ g:
=st^{{1\over p}}st^pst^{{1\over p}}\rangle\cr}$$

\item{(c)} Write $N=2^em$ where $m$ is odd. Let $d$ be any integer
satisfying the congruences $d\equiv 1$ (mod $2^e$), $d\equiv 0$ (mod $m$).
Write $d_2$ and $d_3$ for the multiplicative inverses (mod $N$) of
$2-d$ and $2d+1$. Then
$$\eqalign{{\rm SL}_2(N)=\langle s,t\,|\,&t^N=s^4
=[t^{2^e},st^ms^{-1}]=[g_*,t]=1,\,(st^{-1})^3=s^2,\,
 g_*s=sg_*^{-1},\cr&g_2s=sg_2^{-1},
g_2t=t^{4-3d}g_2, g_3s=sg_3^{-1},g_3t=t^{8d+1}g_3,\ {\rm where}
\cr&g_*:=(st^{1-2d})^3, g_2:=st^{d_2}st^{2-d}st^{d_2},
 g_3:=st^{d_3}st^{2d+1}st^{d_3}\rangle\cr}$$}

Of course $d$ in part (c) is guaranteed to exist, by the Chinese Remainder Theorem.
By e.g.\ `$[g_*,t]=1$' we mean that $g_*$ and $t$ commute.
We have in mind here that $t=\left(\matrix{1&1\cr 0&1}\right)$, $s=\left(
\matrix{0&1\cr -1&0}\right)$, and $g=\left(\matrix{p&0\cr 0&{1\over p}}
\right)\longleftrightarrow G_p$. The monomial
 $st^ms^{-1}$ in (c) will then be $\left(\matrix{1&0\cr -m&1}\right)$. 
Note in (c) that $d_2={1\over 2}(1+d+iN)$ where $d\equiv
1+iN$ (mod $2^{e+1}$), and that $d_3= 1-{2\over 3}(d-jN)$
where $d\equiv jN$ (mod $3^{f+1}$) and $3^f$ is the power of 3
exactly dividing $N$. Since in all these cases SL$_2(N)$
clearly satisfies the given relations, all we must prove here is that we
have included enough relations. We won't use (b) in what follows. 
Our proof of the third presentation exploits the fact that SL$_2(LM)\cong{\rm 
SL}_2(L)\times{\rm SL}_2(M)$ whenever $L$ and $M$ are coprime --- we could
have also used powers of 3 rather than 2. An important
feature of our presentations is that the $g$'s
 correspond to certain automorphisms of
the group SL$_2(N)$,  as we will discuss in \S3 and exploit shortly.

\medskip{\it Proof of the Lemma.} Writing $c:=s^2$, we get in all
three cases that $c$ commutes with $s$ and $st^{-1}$, hence with
everything, and also that $(st)^3=(ts)^3=1$.

To get a presentation for SL$_2(N)$ when $N$ is coprime to some prime
$p$, it is enough to adjoin the relation `$t^N=1$' to any presentation
for the infinite group SL$_2(\Z[{1\over p}])$ where $\Z[{1\over p}]$
denotes the ring $\{{\ell\over p^i}\,|\,\ell,i\in\Z\}$. This important
fact is a quick corollary of a deep theorem by Mennicke [M] showing
that any finite subgroup of SL$_2(\Z[{1\over p}])$ contains a
congruence subgroup; 
the short proof of that corollary is given on p.1433 of [BM]. 
Presentations for SL$_2(\Z[{1\over p}])$ are
given in [M] and most effectively [H2]. In particular, for $p=2$ or 3 the
presentation given in Theorem 5 of [H2] implies that whenever $N$ is
coprime with $p$, SL$_2(N)$ is generated by $x,y$ satisfying 
$$\eqalignno{xy^{-p}x=&\,y^{-p}xy^{-p}&(4a)\cr (xy^{-p}x)^2=&\,(x^py^{-1}x^p)^2
&(4b)\cr x^py^{-1}x^p=&\,y^{-1}x^py^{-1}&(4c)\cr
(xy^{-p}x)^4=&\,x^N=1&(4d) \cr}$$

Put $x=t$ and $y=st^{-{1\over p}}s^{-1}$. It is enough to show, using
the relations given in part (a), that this substitution satisfies the
five relations in (4).

Now, $xy^{-p}x=tsts^{-1}t=s$ and $x^py^{-1}x^p=t^pst^{{1\over
p}}s^{-1}t^p= t^pgt^{-{1\over p}}s=gs$, so equations (4) say
$s=sts^{-1}tsts^{-1}$, $s^2=gsgs$, $gs=st^{{1\over
p}}s^{-1}t^pst^{{1\over p}}s^{-1}$ and $s^4=t^N=1$, all of which
clearly follow from the relations in part (a).

The proof of (b) is similar, and uses the presentation of
SL$_2(\Z[{1\over p}])$ given on p.944 of [H2].

Finally, turn to the most difficult case: part (c).
 Define $S_e=t^dst^dst^dc$, $T_e=t^d$, $S_o=sS_e^{-1}$ and
$T_o=t^{1-d}$, where $c=s^2$. We first want to show $S_e$ and $T_e$ 
commute with both $S_o$ and $T_o$. 

Begin with the observation that the relation $[t^{2^e},st^ms^{-1}]=1$ means,
taking appropriate powers, that $t^{d-1}$ and $st^ds^{-1}=st^dsc$ commute.
Hence $st^dst^dst^d=t^{d-1}st^dstst^d=t^{d-1}st^dt^{-1}sct^{-1}t^d=
t^{d-1}st^{d-1}st^{d-1}c$, and by similar reasoning $t^dst^dst^ds=t^dstst^dst^{d-1}
=t^{d-1}st^{d-1}st^{d-1}c$. Thus $s$ and $S_e$ commute, so so do $S_e$ and $S_o$.
This calculation says $S_o=(st^d)^{-3}=(t^ds)^{-3}=t^{1-d}st^{1-d}st^{1-d}
c$,  hence $T_e(st^d)^3=(t^ds)^3T_e$
and so we get that $S_o$ and $T_e$ also commute.

By the same reasoning, $S_eT_o=t^dst^dstc=tst^dst^{d}c=T_oS_e$,
so $S_e$ and $T_o$ commute. Trivially, $T_e$ and $T_o$ commute. Thus
the subgroups ${\cal G}_e:=\langle S_e,T_e\rangle$ and ${\cal G}_o:=\langle S_o
,T_o\rangle$ commute, and ${\cal G}_e\times {\cal G}_o$ equals the full group 
generated by $s,t$. 
We will be done if we can show ${\cal G}_e$ and ${\cal G}_o$ obey the 
 SL$_2(2^e)$ and SL$_2(m)$ relations, obtained from (a).

Note that $g^2_*=g_*(st^{1-2d})^3=st^{1-2d}g_*^{-1}(st^{1-2d})^2st^{1-2d}t^{2d-1}
s^{-1}=st^{1-2d}t^{2d-1}s^{-1}=1$, so $g_*$ has order 2 and commutes with 
both $s,t$. Define now $\alpha(s)=g_*s$ and
$\alpha(t)=t^{1-2d}$. Because $\alpha(s)$ and $\alpha(t)$ obey all our
relations in (c), $\alpha$ extends to a well-defined group
endomorphism of $\langle s,t\rangle$ --- in fact a group automorphism
since $\alpha^2=id.$
(of course $(1-2d)^2\equiv 1$ (mod $N$)).
Now, hit $t^{1-d}st^{1-d}st^{1-d}c=S_o=(st^d)^{-3}$ with  
$\alpha$: we get $S_o=(g_*st^{-d})^{-3}=(t^ds)^3g_*c=S_o^{-1}g_*c$. 
Thus $S_o^2=g_*c$, and hence $S_o^4=1$. Together with $s^4=1$, we get 
$S_e^4=1$. Now, $(S_eT_e)^3=(S_o^{-1}sT_e)^3=S_o^{-3}(st^d)^3=S_o^{-4}
=1$ and hence also $(S_oT_o)^3=1$.

Finally, $g_2$ and $g_3$ obey the relations for $g$ appearing in 
the SL$_2(m)$ and SL$_2(2^e)$ presentations obtained from (a); that we
have  both $g_2\in{\cal G}_o$ 
and $g_3\in{\cal G}_e$, can be seen by using two more automorphism arguments.
Namely, define $\alpha_2(s)=g_2s$, $\alpha_2(t)=t^{2-d}$, then $\alpha_2$
is seen to define an  automorphism for our group $\langle s,t\rangle$.
Evaluating $\alpha(1)=\alpha_2(S_oT_o)^3$ gives
$1=S_og_2^{-1}T_o^2g_2S_oT_o^2S_og_2^{-1}T_o=S_oT_o^{1\over 2}S_oT^{2}S_o
T^{1\over 2}g_2^{-1}$,
hence $g_2=S_oT_o^{1\over 2}S_oT_o^2S_oT_o^{1\over 2}$ as it should. The 
proof for $g_3$ is identical. \quad \qed\medskip

The best presentation for $N=m$ odd is [CR]:
$${\rm SL}_2(m)=\langle x,y\,|\,x^2=(xy)^3,\,(xy^4xy^{(m+1)/2})^2y^mx^{2k}=1 
\rangle$$
where $k:=[m/3]$ (rounded down). This is optimal, in the sense that at
least 2 generators are needed (since SL$_2(m)$ isn't cyclic) and at
least as many relations as generators are needed (since SL$_2(m)$ is
finite). More generally, for any finite group $G$, the number of
relations minus the number of generators in any presentation must 
at least equal the rank of the `Schur multiplier' $M(G)$ of
$G$. $M(G)$ is always a finite abelian group; for $G={\rm SL}_2(N)$ it
was computed in [B1] and equals $\Z/2\Z$ whenever 4 divides $N$,
otherwise it's trivial. Hence when 4 divides $N$, the best possible is
3 relations, but it seems 5 (the number we use) is the best that has
been achieved thus far in the literature. Of course it goes without saying that the
usefulness of a presentation is not merely determined by the number of
generators and relations.

\smallskip{\it 2.3. The congruence property for general RCFT.}
The following theorem is one of two main results in this paper. 

\medskip {\smcap Theorem 2.} {\it Consider any RCFT. 
 Choose any integer $N$ so that $T^N=I$. Then
our RCFT has the $\Gamma(N)$ congruence property, provided either:\smallskip

\item{(a)} for  $N$ coprime to either $p=2$ or $p=3$,
$G_pT=T^{p^2}G_p$; \smallskip

\item{(b)} for arbitrary $N=2^em$ where $m$ is odd (let $d$ be as in
the Lemma), the following four relations all hold:

\item\item{(i)} $T^{2^e}$ commutes with $ST^mS^{-1}$;

\item\item{(ii)} $G_{2d-1}\,T=T\,G_{2d-1}$;

\item\item{(iii)} $G_{2-d}\,T=T^{4-3d}\,G_{2-d}$; and

\item\item{(iv)} $G_{1+2d}\,T=T^{1+8d}\,G_{1+2d}$}.\medskip

{\it Proof.}  Consider first $N$ in (a). 
Consider the equation $(ST)^3=I$, and apply the Galois automorphism $\si_p$ 
to it. We get $SG_p^{-1}T^pG_pST^pSG_p^{-1}T^p=I$, which we can simplify using
(a) to get $S\,T^{{1\over p}}\,S\,T^p\,S\,T^{{1\over p}}=G_p$.

We find then that the assignments $s\mapsto S$, $t\mapsto T$, $g\mapsto G_p$ 
obey all
the relations in the Lemma, proving Theorem 2(a).

The proof of Theorem 2(b) is  similar.\quad \qed\medskip

In part (b), the matrix
$ST^mS^{-1}$  corresponds to the $2\times 2$ matrix $\left(\matrix{
1&0\cr -m&1}\right)$. We see that $G_{1-2d}$ lies in the centre.

The conditions (a) and (b)(ii),(iii),(iv)  appearing in Theorem 2 are 
surprisingly simple, and we show in Theorem 4 and
\S4 how easy they are to verify in practice. In (i), the roles of $m$ and
$2^e$ can be interchanged if it is more convenient. Also, we could
just as easily use the factorisation $N=3^f\ell$ for gcd$(\ell,3)=1$
and SL$_2(N)\cong {\rm SL}_2(3^f)\times{\rm SL}_2(\ell)$, and make the
appropriate changes to (i)--(iv). 

Consider (i): it is equivalent to the statement
$${\cal U}_{ac}:=\sum_b S_{ab}  S_{c b}^*  
T_{bb}^m=0\ {\rm unless}\ T_{aa}^{2^e}
=T_{cc}^{2^e}\ .\eqno(5)$$
This matrix ${\cal U}$ is symmetric and 
 unitary, and has order $2^e$. Note
that for any $\si_\ell\in{\rm Gal}(\Q[S,T]/\Q)$, $\si_\ell{\cal U}_{ac}
=\epsilon_\ell(a)\,\epsilon_\ell(c)\,({\cal U}^\ell)_{\si a,\si
c}$. To see how (5) can be proved in practice, see section 4. We
will have more to say about ${\cal U}$ shortly.

The converse of Theorem 2 can also be expected to hold: we expect (6a)
below to hold for all $\ell$, and condition (i) is a consequence of the
factorisation SL$_2(N)\cong{\rm SL}_2(m)\times{\rm SL}_2(2^e)$.
In fact, if the characters ch$_a$ were all linearly independent, then by the
Galois action argument of \S3, Theorem 2 would be an `if and only if'.

\smallskip{\it 2.4. Galois and T}.
Suppose now we have an RCFT which may or may not have the congruence property.
In practice (see \S4) it is easy to verify any condition of the form
$$G_\ell T=T^{\ell^2}G_\ell\eqno(6a)$$
or equivalently
$$T_{\si_\ell a,\si_\ell a}=T^{\ell^2}_{aa}\ .\eqno(6b)$$
Let us derive some easy consequences of (6a). Clearly, if $\ell$ obeys (6a), so
does any $\pm\ell^j$.

 Let $M $ be an integer such that the cyclotomic field 
$\Q[\xi_M]$      contains  all entries of 
$S$ and $T$ --- we know it exists by [CG]. Then for any $\ell\in \Z_M^* $
obeying (6a), hitting $(ST)^3=I$ with $\si_\ell$ gives
$$G_\ell=S\,T^{{1\over \ell}}\,S\,T^\ell\,S\,T^{{1\over \ell}}=T^\ell\,S\,
T^{{1\over \ell}}S\,T^\ell\,S\ .\eqno(6c)$$
(6c) has two immediate consequences. Firstly, if $\ell\in\Z_M^*$ obeys
(6c), then the group $\langle S,T\rangle$ has the automorphism
defined by $S\mapsto G_\ell S$, $T\mapsto T^\ell$. 
The reason is that any relation between $S$ and $T$ (i.e.\ monomial $S^aT^b\cdots
S^yT^z=I$) will also be obeyed by $\si_\ell(S)=G_\ell S$ and $\si_\ell T=
T^\ell$ ((6c) is needed to tell us that $\si_\ell(S)\in\langle S,T\rangle$).
Secondly, if the RCFT has the $\Gamma(M)$ congruence property and $\ell\in
\Z_M^*$, then  $G_\ell$ corresponds via $\rho$ to the matrix 
$\left(\matrix{\ell&0\cr 0&{1\over \ell}\cr}\right)\in{\rm SL}_2(M)$.

The following consequences of (6) are valid whether or not the congruence
property holds.

\medskip{\smcap Proposition 3.} {\it Suppose (6b) is valid for all $\ell\in\Z_M^*$.
Let $N$ be the order of $T$: $T^N=I$.

\item{(a)} Then $S_{ab}\in\Q(\xi_N)$ (i.e.\ we may take $M=N$ above).

\item{(b)} Suppose that all `quantum-dimensions'
${S_{a0}\over S_{00}}$ are rational. Then the central charge $c$
is an integer.

\item{(c)} Fix any $b\in\Phi$. Choose any positive integer $K_b$ for which
 all ratios ${S_{ab}\over S_{0b}}$, as `$a$' varies over $\Phi$,
  lie in the cyclotomic field $\Q[\xi_{K_b}]$. Let $M_b$ be
least common multiple of $K_b$ with the order of the root of unity
$T_{bb}$. Then the ratio 
$M_b/K_b$ is a divisor of 24 which is coprime to $K_b$. 
}\medskip

{\it Proof.} If $\si_\ell T= T$, then $T^\ell=T$ so by (6c) we get $G_\ell=I$,
i.e.\ $\si_\ell S=S$, and (a) holds.

Fix any $b\in\Phi$, and define ${\Bbb K}_b$ to be the field generated over
$\Q$ by all ratios ${S_{ab}\over S_{0b}}$ $\forall a\in\Phi$. So
${\Bbb K}_b\subseteq \Q[\xi_{K_b}]$. Of course
${S_{ab}^*\over S_{0b}}={S_{Ca,b}\over S_{0b}}\in{\Bbb K}_b$ and ${1\over
S_{0b}^2}=\sum_a{S_{ab}\,S_{ab}^*\over S_{0b}^2}\in{\Bbb K}_b$. Choose any
$\si=\si_\ell\in{\rm Gal}(\Q(\xi_N)/{\Bbb K}_b)$ --- for instance any
$\ell\equiv 1$ (mod $K_b$) will work. Then $S_{0b}^2=\si(S_{0b})^2=
S_{0,\si b}^2$, so $S_{0b}=s\,S_{0,\si b}$ for some sign $s\in\{\pm 1\}$. 
Hence ${S_{ab}\over S_{0b}}=\si{S_{ab}\over S_{0b}}
={S_{a,\si b}\over S_{0,\si b}}=s\,{S_{a,\si b}\over S_{0b}}$, i.e.\ $S_{ab}
=s\,S_{a,\si b}$ $\forall a$. Unitarity of $S$ forces $b=\si b$ (and $s=+1$). 
Now (6b) gives $T_{bb}^{\ell^2}=T_{bb}$.

In (b) we restrict to $b=0$, and we find that $T^{\ell^2}_{00}=T_{00}$ for all
$\ell\in\Z_N^*$. But $T_{00}=\exp[-2\pi\i\,c/24]$. Let $n$ be the denominator
of the rational number $c/24$. Then $\ell^2\equiv 1$ (mod $n$). But the
`definition of 24' says that that congruence can be satisfied for all
$\ell\in\Z_n^*$, iff $n$ divides 24. Hence $c\in\Z$.

The more general (c) is only slightly more complicated. Consider
$\ell$ coprime to $N$, of the form $\ell=1+K_b\ell'$ for some $\ell'$.
We must have $\ell^2\equiv 1$ (mod $N_b$), where $N_b$ is the order of
the root of unity $T_{bb}$. This means $M_b/K_b$ must divide
$(2+\ell')\ell'$. We are free to choose any $\ell'$ for which
$1+K_b\ell'$ is coprime to $N$; the reader can quickly use that
freedom to prove (c). \qquad\qed\medskip

The `definition of 24' is an easy-to-prove little fact about the
number 24, which explains why that number appears in so many
places. It states that for any $\ell\in\Z$, the relation
$\ell^2\equiv 1$ (mod $n$) holds for all $\ell$ coprime to $n$, iff
$n$ divides 24.

Examples of 3(b) are provided e.g.\ by the finite group orbifolds of holomorphic
theories (see [DVVV], or \S4.3 below), or by the WZW theories
associated to e.g.\ $D_\ell^{(1)}$ level 2 when the rank $\ell$ is a
perfect square.

Equations (6) also have consequences for (5). If (6a) holds for some $\ell$,
then we get the alternate expression $\si_\ell{\cal U}_{ac}=({\cal U}^{1\over
 \ell})_{ac}$. Now, suppose (6a) holds for all Galois automorphisms $\si_\ell$
with $\ell\equiv 1$ (mod $2^e$).
Then $\si_\ell{\cal U}_{ac}={\cal U}_{ac}$ for these $\ell$,
and hence ${\cal U}_{ac}\in\Q[\xi_{2^e}]$.

Let ${\frak C}$ denote the {\it commutant} of the RCFT, i.e.\ the set of all
complex matrices commuting with both $S$ and $T$. This vector space is interesting because
it contains  the
coefficient matrix of the genus-one partition function of the theory. Now, [BI]
discovered the curious fact that the commutant of the affine algebra 
$A_\ell^{(1)}$ level $k$
has a basis $M_1,\ldots,M_n$ consisting of integral matrices. This was
later extended to all affine algebras. We will show
in the next section that this is in fact a generic feature of RCFTs. For any
Galois automorphism $\si$, applying
$\si$ to $M_a=SM_aS^*$ gives that $M_a$ commutes with $G_\si$, and hence
any $M\in{\frak C}$ commutes with all $G_\si$. By the double-commutant
theorem, we know then that each $G_\si$ can be written as a polynomial in 
$S,T$. %In addition, if we knew that in fact $G_\si$ was in the group generated
%by $S$ and $T$, then we would know we could write $G_\si=S^aT^bS^cT^dS^eT^f$
%for some $a,b,\ldots,f\in\Z$ (this property is true for any element in any
%subgroup of $\Gamma$ for which $T$ has finite order).
Equation (6c) finds that polynomial explicitly for us, provided (6a) holds
for that $\si$.

The argument giving (6c) is quite general. In particular, consider now 
 any RCFT --- (6a) and the congruence property may or may not be
satisfied. Write $T_{(\ell)}:=G_\ell T^{{1\over \ell^{2}}}G_\ell^{-1}$; it is 
diagonal, with entries $T_{(\ell)aa}=T_{\si_\ell a,\si_\ell a}^{{1\over
\ell^{2}}}$.
Equations (6) hold iff $T_{(\ell)}=T$. The argument giving (6c) now yields
$$G_\ell=ST^{{1\over \ell}}ST_{(\ell)}^\ell ST^{{1\over \ell}}=T^\ell S
T_{({1\over \ell})}^{{1\over \ell}}ST^\ell S=ST_{({1\over \ell})}^{{1\over \ell}}
ST^\ell ST_{({1\over \ell})}^{{1\over \ell}}
=T_{(\ell)}^\ell ST^{{1\over \ell}}ST_{(\ell)}^\ell S\eqno(7)$$
Again, (7) holds in complete generality. Since by construction the
matrix  $G_\ell$ is real, we also get $G_\ell=ST^{*{1\over
\ell}}ST^{*\ell}_{(\ell)} ST^{*{1\over \ell}}C$ etc. We are now ready for our
second main result.

\smallskip{\it 2.5. The congruence property for RCFTs when $T$ has odd
order}.

\medskip{\smcap Theorem 4.} {\it Consider any RCFT. Let $N$ be the order of $T$:
i.e.\ $T^N=I$. If $N$ is odd, then the RCFT obeys the $\Gamma(N)$ congruence property.}

\medskip{\it Proof.} It was shown in [By] that for any $a,b\in\Phi$,
the number
$${\cal Z}(a,b):=T_{00}^{{1\over 2}}T_{bb}^{*{1\over 2}}\sum_{x,y\in\Phi}
N_{xy}^aS_{bx}S_{0y}T^2_{yy}T^{-2}_{xx}\eqno(8)$$
is an integer, among other things. We can interpret all fractions as integers
(mod $N$), as usual. Consider the matrix ${\cal Z}^{(a)}:=
T^{{1\over 2}}ST^2N_aT^{*2}ST^{*{1\over 2}}$ where $N_a$ is the fusion
matrix $(N_a)_{bc}=N_{ab}^c$. Write $T_h:=T_{({1\over 2})}$; we know from (7) that 
$G_2=ST^{*{1\over 2}}_hST^{*2}ST^{*{1\over 2}}_hC$ and $G_{{1\over 2}}=T^{{1\over 2}
}_hST^2ST^{{1\over 2}}_hS$. Thus we can write
$$\eqalignno{{\cal Z}^{(a)}=&\,T^{{1\over 2}}T_h^{*{1\over 2}}(G_{{1\over 2}}
ST^{*{1\over 2}}_h S)N_a(ST_h^{{1\over 2}}SG_2C)T^{{1\over 2}}_hT^{*{1\over
2}}&\cr&=T^{{1\over 2}}T_h^{*{1\over 2}}G_{{1\over
2}}ST^{*{1\over 2}}_hD_aT^{{1\over 2}}_hSG_2T^{{1\over 2}}_hT^{*{1\over 2}}&}$$
where $D_a$ is the diagonal matrix with entries $(D_a)_{xx}={S_{ax}\over S_{0x}}$
(we used the Verlinde formula (1) here). Since $D_a$ and $T_h$ are both
diagonal, they commute and we get
$${\cal Z}^{(a)}=T^{{1\over 2}}T^{*{1\over 2}}_hG_{{1\over 2}}N_{Ca}G_2
T_h^{{1\over 2}}T^{*{1\over 2}}\ .$$
Put $\alpha=T_{00}T_{h\,00}^*$, then (8) tells us 
$$\sqrt{\alpha}\,
\epsilon_{{1\over 2}}\!(0)\,\epsilon_{{1\over
2}}\!(b)\,N_{Ca,\si_{{1\over 2}} \!0}^{
\si_{{1\over 2}}b}\,T^{{1\over 2}}_{h\,bb}\,T^{*{1\over 2}}_{bb}\in \Z$$
for all $a,b\in\Phi$. Hence $\alpha \,T_{h\,bb}\,T_{bb}^*\in\Q$ for all $b$.
But it is also an $N$th root of unity, and $N$ is odd, so we find $T=\alpha T_h
=\alpha T_{({1\over 2})}$.
Therefore $$\eqalign{I&=G_2G_{{1\over 2}}=(\alpha^{-1}ST^{{1\over 2}}ST^2S
T^{{1\over 2}})(
\alpha^{-1}T^{{1\over 2}}ST^2ST^{{1\over 2}}S)\cr&=\alpha^{-2}ST^{{1\over 2}}
ST^2(T^{-1}ST^{-1}
C)T^2ST^{{1\over 2}}S=\alpha^{-2}ST^{{1\over 2}}(T^{-1})T^{{1\over 2}}SC=
\alpha^{-2}I}$$
where we made repeated use of $(ST)^3=I$. Hence $\alpha=1$ and $T=T_{({1\over 2})}$,
and by Theorem 2(a) we are done.\qquad \qed\medskip

For clarity, let us repeat the properties of RCFTs used in deriving Theorems 2
 and 4. We used the facts that $S$ and $T$ are both unitary and
 symmetric, and that $T$ is
diagonal and of finite order. We used Verlinde's formula (1), or at 
least the fact that that formula always yields rational numbers (it should
in fact give nonnegative integers). Finally, the derivation of (8) in [By]
required elementary properties of the ($N_{aa}^b\times N_{aa}^b$) braiding
matrix ${\cal R}_{aa}^{[b]}$ ($=\Omega_{aa}^b$ in the notation of [MS]): he 
needs its square and trace.
All of these are standard properties of RCFTs [MS].

The condition that $T$ has odd order, can be rephrased in terms of
more elementary quantities of RCFT, namely the central charge $c$ and
the conformal weights $h_i$. First note that we can write any
rational number $r$ uniquely as a product of prime numbers raised to certain
integral powers: $r=\prod p_i^{a_i}$. Let $t(r)$ denote the `two-ness'
of $r$, i.e.\ the exponent of 2 in this prime decomposition of $r$. For
example, $t(2.4)=2, t({5\over 3})=0, t(33.5)=-1$. Then $T$ has odd
order, iff $t(c)\ge 3$ and $t(h_i)\ge 0$ for all $i$.

\smallskip{\it 2.6. Further problems}.
Additional directions suggested by this work are the following. Roughly,
they involve the interplay between three topics: the congruence property, 
equations (6), and the number fields $\Q[S]$ and $\Q[T]$. See also Conjecture
3 next section.

\smallskip\item{$\bullet$} Can we prove (6a) uniformly for all RCFT?
Take any Galois automorphism $\si_\ell$, and apply $\si_\ell^2$
to the equation $(ST)^3=I$. We get $(G_\ell S G_\ell^{-1}T^{\ell^2})^3=I$,
i.e.\
$$(ST')^3=I\ ,\ {\rm where}\ T':=G_\ell^{-1}T^{\ell^2}G_\ell=T_{({1\over \ell})}
\eqno(9)$$
We want to show $T'=T$. The point is that the equation $(ST)^3=I$ is strong
and almost uniquely determines $T$ from $S$, given the facts that $T$ is diagonal
and of finite order. It could be useful to understand to what extent
these conditions do determine $T$. Clearly,
given one solution $T$ of $(ST)^3=I$, we can always 
obtain another by multiplying the matrix by a third root of unity. If $S$ is real, then
$T^*$ will be another solution. For a given $S$ there can be additional
`sporadic' solutions for $T$, but there doesn't seem to be that many. Thus (9)
can be interpreted as suggesting that (6a) will generically hold.

\smallskip\item{$\bullet$} 
In the case (b) of Theorem 2 (i.e.\ $N$ a multiple of 6),
if we drop e.g.\ condition (i), can we  at least show that the matrices
$S$ and $T$ generate a {\it finite} group? 

\smallskip\item{$\bullet$} When  we have the
congruence property, we should be able to prove equations (6) 
for all $\ell$. An approach is discussed next section, but it requires 
additional knowledge of the RCFT characters.

\smallskip\item{$\bullet$} Conversely, it is tempting to believe that if 
equations (6) hold for all $\ell$, then we should have the congruence 
property. Of course the only question is proving (5) when $N$ is
a multiple of 6. We could accomplish this if we knew that for any $N$,
SL$_2(N)$ has a presentation of the form
$$\eqalign{{\rm SL}_2(N)\cong\langle s,t\,|\,&t^N=s^4=1,(st^{-1})^3=s^2,\,
g_as=sg_a^{-1},\, g_at=t^{a^2}g_a,\, g_ag_b=g_bg_a\ \cr
&{\rm where}\ g_a:=st^{{1\over a}}st^ast^{{1\over
a}}\ {\rm for\ all}\ a,b\in\Z_N^* \rangle}$$
Lemma 1 says that this convenient but terribly inefficient presentation works
whenever $N$ is not a multiple of 6. We  expect this
presentation to hold for all $N$, but we don't know a proof. Can
anyone help us?

\smallskip\item{$\bullet$} To our knowledge, the powers of the matrix
$STS^{-1}=\rho(\left(\matrix{1&0\cr -1&1\cr}\right))$ have been largely
ignored in RCFT and the studies of the associated representations of the 
modular group,
but they do play a role in the mathematical theory, as Theorem 2(b)
demonstrates. We give an example  of how to study this, at the
beginning of section 4.

\smallskip\item{$\bullet$} We expect that the $\Gamma(N)$ congruence property 
should imply that all entries $S_{ab}$ lie in the cyclotomic field $\Q[\xi_N]$.
We know [CG] the $S_{ab}$ lie in {\it some} cyclotomic field, but the
question here is whether they lie in that specific one.
We know this will hold if (6) holds for all $\ell$. We know from basic
results in modular functions theory (see e.g.\ [L]) that the Fourier coefficients
for each ch$_a(-1/\tau)$ will lie in $\Q[\xi_N]$ provided ch$_a$ is fixed
by $\Gamma(N)$. This last statement is sometimes (see e.g.\ [ES]) taken 
(prematurely it seems to us)
to imply that each $S_{ab}\in\Q[\xi_N]$ when the $\Gamma(N)$ congruence 
property holds. The problem is the usual one: the characters ch$_a(\tau)$
in an RCFT won't in general be linearly independent, so we can't read off
$S$ from (2a). We will return to this question next section.

\smallskip\item{$\bullet$} The automorphism group of each PSL$_2(p^n)$
is known [MD], and from this it is straightforward to obtain the
automorphism group of each SL$_2(p^n)$. In particular, for $p>5$ the
automorphisms are generated by the Galois ones $s\mapsto g_\ell s$,
$t\mapsto t^\ell$ for each $\ell$ coprime to $p$, together with the
inner ones $a\mapsto bab^{-1}$ for each $b\in{\rm SL}_2(p^n)$ --- in
particular the outer automorphism group Out$({\rm SL}_2(p^n))\cong
\Z/2\Z$.  There
are additional automorphisms for $p=2,3,5$. Perhaps it can be hoped that for
some $N$, there will be some kind of generalisation of the Galois
action (3) corresponding to nonGalois automorphisms in Out$({\rm
SL}_2(N))$. We briefly return to this next section.

\smallskip\item{$\bullet$} It is clearly desirable to try to extend our
results to any rational vertex operator algebra. The main barrier
 is that Verlinde's formula (1) is not yet
known to give rationals there (or even to be defined!). Theorem 1 in
[DLM]  shows that for any
rational VOA obeying the technical finiteness condition `$C_2$', $T$ will
have finite order. Again subject to the $C_2$ condition, Theorem 5.3.3 in
[Z] (see also Theorem 3 in [DLM]) shows that the characters ch$_a(\tau)$
for a rational VOA will all be holomorphic in the upper half-plane, and
define a representation of SL$_2(\Z)$. This $C_2$ condition is conjectured
to hold for all rational VOAs. A discussion of VOAs in RCFT is provided
e.g.\ by [H3].

\bigskip{{\bf 3. Explaining some curiousities}}\medskip

See also [B] for some related comments.

\smallskip{\it 3.1. The RCFT Galois action reinterpreted.} 
 The Galois action [CG] in RCFT seems somewhat mysterious, but actually
it is a special case of one known since early this century.
A cyclotomic Galois action arises naturally in modular functions (see
e.g.\ Chapter 6 of [L]).
In particular, let
 $f(\tau)=  q^c  \   \sum_{n=0}^\infty a_nq^{n/N }$ be a modular 
function for $\Gamma(N)$, with
coefficients $a_n$ in $\Q[\xi_N]$. Choose any $\ell\in\Z_{N}^* $ and write
$h_\ell$ for any matrix in $\Gamma$ congruent (mod $N$) to 
$\left(\matrix{\ell&0\cr 0&{1\over \ell}}\right)$. Then we get
the remarkable formula [L]
$$f(h_\ell\tau)=  q^c  \  
 \sum_{n=0}^{\infty}\si_\ell(a_n)\,q^{n/N }\eqno(10a)$$
which will also be a modular function for $\Gamma(N)$.

Now let ch$_a(\tau)$ be the characters for an RCFT obeying the congruence
property. Apply this Galois action to ch$_a(s\tau)$: we find
$${\rm ch}_a(h_\ell s\tau)=\sum_{b\in\Phi}\si_\ell(S_{ab})\,{\rm ch}_b(\tau)\eqno(10b)$$
and hence $\rho(h_\ell)=G_\ell$, as in (6c).

As mentioned in \S2, arguments of these kind break down in most RCFTs, 
because they assume that the RCFT characters are linearly independent. Introducing
additional variables into VOA (hence RCFT) characters can be done quite
generally, as will be discussed more fully in [GG]. Variants of Jacobi
forms [EZ] often arise in this way (these have a 
`linear' $\vec{z}$-dependence
as well as the `quadratic' $\tau$ dependence, in analogy with the function
theta $\theta(\tau,z)$ discussed in \S2). Jacobi forms behave essentially
the same as modular functions, and the analogue of the cyclotomic Galois action was
worked out in [Be] (he restricted to a single variable $z$, but his argument 
extends to vectors $\vec{z}$). We find that the above $\rho(h_\ell)=G_\ell$
observation carries over. Also, the representation of $\Gamma$ (or
SL$_2(N)$) arising from these Jacobi functions will have matrix entries
in the field $\Q[\xi_N]$. 

SL$_2(N)$ has automorphisms defined by $s\mapsto h_\ell s$, $t\mapsto t^\ell$,
for any $\ell\in\Z_{N}^* $. Looking at $t$, it is easy to show that the
$\ell$th automorphism is {\it inner} (i.e.\ given by $M\mapsto VMV^{-1}$
for some $V\in{\rm SL}_2(N)$), iff $\ell$ is a perfect square mod $N$.
These automorphisms are precisely the Galois actions on modular functions
for congruent subgroups, and are precisely the Galois actions in RCFT
when that RCFT obeys the congruence property. Thus in hindsight it seems
that we should have taken the presence of the RCFT Galois action of [CG]
as a strong hint that the $\Gamma$ representation $\rho$ there factors through
a congruence subgroup, or, if we already suspected that $\rho$
factors through one, then we should have started looking
for an action on our RCFT data of the Galois group $\Z_{N}^* $.

This discussion, and the results of our paper, lead us to propose the
following strengthening of Conjecture 1:\medskip

{\smcap Conjecture 3}. {\it All RCFTs have the $\Gamma(N)$ congruence property,
where $N$ is the order of $T$: in particular $\left(\matrix{0&1\cr -1&0}
\right)\mapsto S$, $\left(\matrix{1&1\cr 0&1}\right)\mapsto T$ defines
a representation of SL$_2(N)$. In addition, each entry $S_{ab}$ lies
in $\Q[\xi_N]$, 
where $\xi_N:=\exp[2\pi\i/N]$, and equations (6) hold for all $\si\in{\rm
Gal}(\Q[\xi_N]/\Q)$.}\medskip

\smallskip{\it 3.2. Integral bases for commutants.}
The congruence property can also be used to provide an explanation
for the existence of an integral basis for the commutant (see the discussion
near the end of \S2.4). In particular:

\medskip {\smcap Proposition 5}. {\it Suppose the matrices $S$ and
$T$ obey both the
$\Gamma(N)$ congruence property as well as equations (6). Let ${\frak C}$ 
denote all the complex matrices
$M$ commuting with both $S$ and $T$. Then ${\frak C}$ has a basis (over ${\Bbb C}$)
consisting of integral matrices.}\medskip

{\it Proof}. We know from Proposition 3 that the entries of $S,T$
lie in the cyclotomic field $\Q[\xi_N]$. Certainly ${\frak C}$ will have a basis
$M_1,\ldots,M_n$ with entries from $\Q[\xi_N]$. Rescaling $M_1$ appropriately,
we may assume some entry $(M_1)_{i_1,j_1}$ equals 1. Subtracting $M_1$ from
the other $M_a$, we can assume all other $(M_a)_{i_1,j_1}=0$. Continue inductively
in this way:
we get indices $(i_1,j_1),\ldots,(i_n,j_n)$ such that $(M_a)_{i_b,j_b}=\delta_{ab}$.

Choose any $M\in{\frak C}$ with entries from $\Q[\xi_N]$, and any 
Galois automorphism $\si=\si_\ell$. Then $M$ commutes with the matrix $G_\si$, by
(6c). From the calculations 
$$\eqalignno{S(\si M)S^{-1}=&\,\si((SG_\si)\, M(G_\si^{-1}S^{-1}))
=\si M\ ,&\cr T(\si_\ell M)T^{-1}=&\,\si_\ell(T^{1\over\ell}MT^{-{1\over \ell}})
=\si_\ell M\ ,&\cr}$$
 we see that $\si M$ also lies in ${\frak C}$. 

Now let $\overline{M}_a=\sum\si M_a$, where the sum is over all $\si\in
{\rm Gal}(\Q[\xi_N]/\Q)$. Then $\overline{M}_a$ will have rational entries,
they will be linearly independent, and they will all lie in ${\frak C}$.
Hence they constitute a basis for ${\frak C}$.\quad \qed\medskip

{\it 3.3. The observations of Bantay.}
We will conclude with a discussion of the remarkable observation in [By]:
the numbers ${\cal Z}(a,b)$ defined in (8) are integers, congruent to
$N_{aa}^b$ (mod 2), and $|{\cal Z}(a,b)|\le N_{aa}^b$. This observation
seems highly nontrivial: for instance taking $a=0$ tells us that
$$\sqrt{T_{00}T_{bb}^*}\,\sum_dS_{db}S_{0d}=\pm\delta_{0b}$$ 
where the sum is
over all self-conjugate $d$ (i.e.\ all $d=Cd$), and this seems difficult to
prove without using either the methods of [By] or this paper.

That the ${\cal Z}(a,b)$ are integers played an important role in our
proof of Theorem 4. Assuming $T$ has odd order, we can turn the Theorem 4
argument around now, and what we find is ${\cal Z}^{(a)}=G_{{1\over 2}}
N_{Ca}G_2$, where the matrix ${\cal Z}^{(a)}$ is defined in the proof
of Theorem 4. ${\cal Z}(a,b)$ is the $(0,Cb)$ entry of ${\cal Z}^{(Ca)}$,
namely
$${\cal Z}(a,b)=\epsilon_{{1\over 2}}(0)\,\epsilon_{{1\over 2}}(b)\,
N_{a,\sigma 0}^{C\sigma b}$$
where for readability 
we write $\sigma$ for the Galois permutation $\sigma_{{1\over 2}}$.
So ${\cal Z}(a,b)$ will indeed always be an integer, at least when $N$ is
odd. This is a consequence of the fact that the $N$ odd RCFT must obey the
congruence property (i.e.\ Theorem 4). The remainder of Bantay's observation now reduces to
$$N_{a,\sigma 0}^{\sigma b}\equiv N_{aa}^{Cb}\quad
({\rm mod}\ 2)\qquad{\rm and}\qquad 
N_{a,\sigma 0}^{\sigma b}\le N_{aa}^{Cb}$$
which cannot be proved using the methods of this paper. Nevertheless,
we see that Bantay's `Frobenius-Schur indicator' ${\cal Z}(a,0)$, which equals
$0$ or $\pm 1$ if $Ca\ne a$ or $a$ is real/pseudo-real, respectively, will never
be negative here. Hence:

\medskip{\smcap Corollary 6.} {\it When $T$ has odd order, the RCFT will have no 
pseudo-real primary fields.}\medskip

In addition, for $N$ odd we get the surprising fact that there exists
an $a\in\Phi$ (namely $a=\sigma_{{1\over 2}}0$) with the property that
$N_{aa}^b\ne 0$ iff $b=Cb$, in which case $N_{aa}^b=1$. In a colourful
phrase suggested to us by M.A.\ Walton, 
\smallskip\centerline{the sum of the self-conjugate primary fields has a
(fusion) square-root!}\smallskip
For a concrete
example, the affine algebra $A_2^{(1)}$ at even level $k$ has $N$ odd,
and $\sigma_{{1\over 2}}0$ there equals the weight $(0,{k\over 2},{k\over 2})$.
Hence we get the fusion rules
$$(0,{k\over 2},{k\over 2})\times (0,{k\over 2},{k\over 2})=
(k,0,0),(k-2,1,1),\ldots,(0,{k\over 2},{k\over 2})\ ,$$
with all  multiplicities equal to 1.

Note that in fact we have shown a little more: {\it any} entry of ${\cal Z}^{(a)}$
is manifestly an integer. Thus this aspect of [By] can be generalised, at least
for odd $N$: the quantities
$${\cal Z}(a,b,d):=\sqrt{T_{dd}T^*_{bb}}\sum_{x,y\in\Phi}N_{xy}^aS_{bx}S_{dy}
T^2_{yy}T^{*2}_{xx}=\epsilon_{{1\over 2}}(d)\,\epsilon_{{1\over 2}}(b)\,
N_{a,\sigma d}^{C\sigma b}$$
will always be integers, for any $a,b,d\in\Phi$ (the case $d=0$ reduces
to [By]). Also, there is nothing special about the number `2' here. In particular,
let $\ell$ be coprime to the order $N$ of $T$ (which we no longer assume to
be odd), and assume that $\ell$ satisfies (6a). Then the numbers 
$${\cal Z}_\ell
(a,b,d):=T^{{1\over \ell}}_{dd}T^{*{1\over\ell}}_{bb}\sum_{x,y\in\Phi}
N_{xy}^aS_{bx}S_{dy}T^\ell_{yy}T^{*\ell}_{xx}$$ will all be integral.
It would be very interesting to find an interpretation for those numbers.

\bigskip {{\bf 4. Examples}}\medskip

{\it 4.1. Lattice theories}.
Consider first the RCFT corresponding to a single compactified boson,
or equivalently a 1-dimensional lattice theory or U(1) theory. Take
the lattice to be $\sqrt{n}\Z$ where $n$ is an even integer. The
primary fields are labelled by $a\in\{0,1,\ldots,n-1\}=\Phi$. The
(full-variable) character for $a\in\Phi$ is proportional to
$\Psi_a(n\tau,\sqrt{n}z)=\sum_{m\in\Z}\exp[2\pi\i\sqrt{n}\,(m+{a\over
n})z +n\pi\i\,(m+{a\over n})^2\tau]$. The transformation law of
the $\Psi_a$ can be read off from that of
$\theta_3(\tau,z)=\Psi_0(\tau,z)$, first found by Poisson and Jacobi
({\it c}.\ 1830). The resulting $S$ and $T$ matrices are
$$S_{ab}={1\over \sqrt{n}}\exp[2\pi\i {ab\over n}]\ ,\qquad
T_{ab}=\exp[\pi\i{a^2\over n}-\pi\i{1\over 12}]\,\delta_{ab}\ .$$
We already know this theory satisfies the congruence property, because
the characters are theta functions. However let's try to see it from
Theorem 2. The order $N$ here is the least-common-multiple of 24 and
$2n$. The Galois permutation is $\sigma_\ell a=\ell a$ taken mod $n$,
for any $\ell$ coprime to $n$. Equation (6b) then is obviously
satisfied, so it suffices to verify (5). Now,
${\cal U}_{ac}={\exp[-\pi\i m/12]\over n}\,S(m,2a-2c,n)$, where
$S(a,b,c)$ is the {\it generalised Gauss sum} 
$$S(a,b,c):=\sum_{k=0}^{c-1}\exp[\pi\i\,(ak^2+bk)/c]\ .$$
$S(a,b,c)$ obeys an important symmetry, called {\it reciprocity}, due
originally to Genocchi (1852) --- see e.g.\ [Bt] for a modern proof
and generalisation. In particular,
$$S(a,b,c)=\sqrt{\left|{c\over a}\right|}\exp[\pi\i\,\{{\rm sgn}(ac)-b^2/ac\}/4]
\,S(-c,-b,a)\ .$$
Applying this to ${\cal U}_{ac}$, we find 
$${\cal U}_{ac}={\exp[\pi\i \,(-m/12+1/4-(a-c)^2/mn)]\over \sqrt{mn}}
\sum_{b=0}^{m-1}e^{2\pi\i\,(c-a)b/m}$$
using the facts that $m$ divides $n$, and $m$ is odd and $n$
even. Hence ${\cal U}_{ac}\ne 0$ iff $m$ divides $a-c$, which implies
that  $m$ divides $a^2-c^2$, and we are done.

Note that this gives an immediate proof that for any $n$-dimensional even Euclidean lattice
$\Lambda$, and any vector $g\in\Lambda^*$ (the dual lattice), the
Jacobi theta function $$\Theta(g+\Lambda)(\tau,z):=
\sum_{x\in g+\Lambda}\exp[\pi\i x^2\tau+2\pi\i x\cdot z]$$
will be fixed by some $\Gamma(N)$ (up to the usual factors). Indeed,
we can express $\Theta(g+\Lambda)$ as a homogeneous degree $n$ polynomial in the 1-dimensional
$\Psi_a$ by using Gram-Schmidt to find in $\Lambda$ an $n$-dimensional
orthogonal sublattice.

\smallskip{\it 4.2. Affine theories}.
Consider next any affine algebra $X_r^{(1)}$  at level $k$. There is
associated  a well-known 
representation $S,T$ of $\Gamma$ [KP] --- the primaries are labelled
by the highest-weights $\la\in P_+^k(X_r)$, and for instance we have
$$T_{\la\la}=\exp[\pi\i\,{(\la+\rho)^2\over k+h^\vee}-\pi\i \,{{\rm dim}
\,X_r\over 12}] \ .$$
Because of the (Jacobi) theta function expression for the affine
characters, we again know this representation of $\Gamma$ factors through 
$\Gamma(N)$,
where $N$ can be taken to be $(k+h^\vee)n$ for some $n$ (e.g.\ $n=24\,(r+1)$ 
works for $A_r^{(1)}$).
 Provided $N$ is not a multiple of 6 (e.g.\ $A_{r}^{(1)}$ level $k$
when $k$ and $r$ are
both even) this also follows from Theorem 2, as we'll now see.

Equations (6) are trivial to verify for the affine
algebras: for any Galois automorphism, 
$\si_\ell(\la)$ can be interpreted as the unique weight $\la^+\in P_+^k(X_r)$
for which 
$$\la^++\rho=w(\ell(\la+\rho))+(k+h^\vee)\alpha$$ where $w$ lies in
the (finite) Weyl group and $\alpha$ lies in the coroot lattice of $X_r$.
Hence the norms $\ell^2(\la+\rho)^2$  and $(\la^++\rho)^2$ are congruent
mod $2(k+h^\vee)$, so relation (6b) holds --- the constant factor $\exp[-\pi\i
\,{\rm dim}\, X_r/12]$ in $T_{\la\la}$ causes no problems,
because of the `definition of 24':
$${\rm gcd}(\ell,24)=1\ \Rightarrow\ \ell^2\equiv 1\ ({\rm mod}\ 24)\ .$$

\smallskip{\it 4.3. Orbifold theories}.
Another important, and in many ways behaviourally opposite, example
of RCFT modular data, is associated to any finite (discrete) group $G$ [DVVV]. 
For any $a\in G$, the conjugacy class associated
to $a$ is the set of all elements of the form $g^{-1}ag$. 
Fix a set $R$ consisting  
of one representative for each conjugacy class of $G$. By $C_G(a)$
we mean the {centraliser} of $a$ in
$G$: i.e. the subgroup consisting of all elements in $G$ which commute with 
$a$. The primary fields here are pairs $(a,\chi)$, where $a\in R$ and 
$\chi$ is the character of an irreducible representation of $C_G(a)$. Here
$$T_{(a,\chi),(b,\chi')}=\delta_{a,b}\,\delta_{\chi,\chi'}{\chi(a)\over
\chi(e)}$$
where  $e$ is the identity
of $G$. This data corresponds to the RCFT obtained by orbifolding a
holomorphic RCFT by $G$. 

The order $N$ of $T$ here is the exponent of $G$ (i.e.\ the smallest positive
integer such that $a^N=e$ for all $a\in G$). Again, it was already
known that the corresponding modular data factors through $\Gamma(N)$.
One way to see this follows from the treatment in [KSSB]. Let ${\cal C}(G)$
be the space of all functions $f:G\times G\rightarrow \C$ satisfying two
conditions:
$f(g,h)=0$ unless $gh=hg$, and $f(aga^{-1},aha^{-1})=f(g,h)$ $\forall a\in G$. 
$\Gamma$ acts
on these by $f((g,h)\left(\matrix{a&b\cr c&d}\right))=f(g^ah^c,g^bh^d)$.
(This strange-looking action of $\Gamma$ on $G\times G$ is related to
the natural action of $\Gamma$ on the fundamental group $\pi_1({\rm
torus})\cong\Z^2$.) This turns out to be the $\Gamma$ representation appearing
 in the RCFT --- this construction shows that  its kernel clearly contains 
$\Gamma(N)$.

Once again, it is immediate that (6b) holds.
 In particular, the Galois permutation (3) takes
$(a,\chi)$ to $(b a^\ell b^{-1},\si_\ell\chi^{b^{-1}})$, where $b\in G$ is chosen
so that $ba^\ell b^{-1}\in R$, and $\chi^g$ denotes
the function $\chi^g(h):=\chi(ghg^{-1})$. Thus $T_{\si(a,\chi),\si(a,\chi)}
=\si^2T_{(a,\chi),(a,\chi)}$, which is (6b).

\smallskip{\it 4.4. The congruence test of Hsu}.
Finally, we should point out that Lemma 1 provides a presentation for
SL$_2(N)$ which may be more natural and hence useful than the one given in [H]
 for $N$ even but not a power of 2.  In this case,
Hsu [H] combined the generators and relations of Mennicke and Behr [M,BM]
to obtain the following presentation of SL$_2(N)$. Write $N=2^em$ where $m$
is odd, and define $d$ as in Lemma 1. Then a presentation for SL$_2(N)$ is
$$\eqalignno{\langle L,R\,|\,& L^N=[a,r]=[b,l]
=(ab^{-1}a)^4=(lr^{-1}l)^4=1, (ab^{-1}a)^2=(b^{-1}a)^3=(b^2a^{-{1\over 2}})^3, 
&\cr& (lr^{-1}l)^2=(r^{-1}l)^3=(sr^5lr^{-1}l)^3, (lr^{-1}l)^{-1}s(lr^{-1}l)=
s^{-1}, s^{-1}rs=r^{25},&\cr&\qquad{\rm where}\
a=L^{1-d}, b=R^{1-d}, l=L^d, r=R^d, s=l^{20}r^{{1\over 5}}l^{-4}r^{-1}\rangle
&\cr}$$
  
The advantage of ours, perhaps, is that our relations involve the
automorphisms of the group and so should be easier to identify and
verify in practice. This may
permit a practical simplification of the congruence subgroup test in [H].
We also have one fewer relation.

Fewer relations even than Lemma 1(c) would be obtained by using the
2 relation presentation of SL$_2({\rm odd})$ in [CR], together with
e.g.\ our 5 relation presentation of SL$_2(2^e)$. It should be pointed
out that provided $N$ is not a multiple of 210, a much better
presentation of SL$_2(N)$ (with at most 6 relations) is given by our
Lemma 1(a),(b).

\bigskip\noindent{{\bf Acknowledgements}}. T.G.\ thanks  C.\  
Cummins, J.\ McKay, D.\ McNeilly, P.\ Moree and M.\ Walton for various
conversations, and A.C.\  thanks D.\  Altschuler, M.\ Bauer, J.\  Lascoux, 
 J.\  Wolfart, and J.-B.\ Zuber. We
 both thank W.\ Nahm and P.\ Ruelle. This paper was written in part  
  while
both authors were visiting IHES, whom we thank
for their generous hospitality. The research of T.G.\ was  supported
in  part by NSERC.

\vfill\eject
%\bigskip\bigskip
\centerline{{\bf Bibliography}}\medskip

\item{[AM]} {\smcap G.\ Anderson and G.\ Moore}, {``Rationality in conformal field
theory''}, Commun.\ Math.\ Phys.\ {\bf 117} (1988), 441--450;

\item{} {\smcap C.\ Vafa}, {``Towards classification of conformal field
theories''}, Phys.\ Lett.\ {\bf B206} (1988), 421--426.

\item{[AS]} {\smcap A.\ O.\ L.\ Atkin and H.\ P.\ F.\ Swinnerton-Dyer}, 
{``Modular forms on noncongruence subgroups''}, In: Proc.\ Symp.\
Pure Math.\ {\bf 19} (AMS, Providence, 1971), ed.\ by T.\ S.\ Motzkin, pp.1--26.

\item{[By]} {\smcap P.\ Bantay}, {``The Frobenius-Schur indicator in conformal
field theory''}, Phys.\ Lett.\ {\bf B394} (1997), 87--88.

\item{[B]} {\smcap M.\ Bauer}, {``Galois actions for genus one rational conformal 
field theories''}, In: {\it The Mathematical Beauty of Physics}, 
(World Scientific, 1997) pp.\ 152--186.

\item{[BCIR]} {\smcap M.\ Bauer, A.\ Coste, C.\ Itzykson, and P.\ Ruelle},
{``Comments on the links
between ${\rm su}(3)$ modular invariants, simple factors in the Jacobian of
 Fermat curves, and rational triangular billiards''}, J.\ Geom.\
Phys.\ {\bf 22} (1997), 134--189.

\item{[BI]} {\smcap M.\ Bauer and C.\ Itzykson}, 
{``Modular transformations of ${\rm SU}(N)$ affine
characters and their commutant''}, Commun.\ Math.\ Phys.\ {\bf 127} (1990),
617--636.

\item{[BM]} {\smcap H.\ Behr and J.\ Mennicke},
{``A presentation of the groups ${\rm PSL}(2,\,p)$''},
 Can.\ J.\ Math.\ {\bf 20} (1968), 1432--1438.

\item{[Bt]} {\smcap B.C.\ Berndt}, {``On Gaussian sums and other exponential
sums with periodic coefficients''}, Duke Math.\ J.\ {\bf 40} (1973), 145--156.

\item{[Be]} {\smcap R.\ Berndt},
{``Zur Arithmetik der elliptischen Funktionenk\" orper h\" oherer
Stufe''}, J.\ Reine Angew.\ Math.\ {\bf 326} (1981), 79--94.

\item{[Bl]} {\smcap F.\ Beyl}, 
{``The Schur multiplicator of ${\rm SL}(2,{Z}/m{Z})$ and the
congruence subgroup property''}, Math.\ Z.\ {\bf 191} (1986), 23--42.

\item{[CR]} {\smcap C.\ M.\ Campbell and E.\ F.\ Robertson}, {``A deficiency
zero presentation of SL(2,p)''}, Bull.\ London Math.\ Soc.\ {\bf 12}
(1980), 17--20.
   
\item{[CG]} {\smcap A.\ Coste and T.\ Gannon},
{``Remarks on Galois symmetry in RCFT''},  Phys.\ Lett.\ {\bf B323} (1994), 
316--321.
        
\item{[C]} {\smcap A.\ Coste},  
{``Investigations sur les caract\`eres de Kac Moody et quelques 
 quotients de  ${\rm SL}(2,{Z} ) $''}, IHES preprint P/97/78, to 
 e-appear, and refs therein. 
      
\item{[CM]} {\smcap H.\ M.\ S.\ Coxeter and W.\ O.\ J.\ Moser}, {\it Generators and Relations
for Discrete Groups}, 4th edn. (Springer, Berlin, 1980).

\item{[DVVV]} {\smcap R.\ Dijkgraaf, C.\ Vafa, E.\ Verlinde, and H.\ Verlinde}, 
{``The operator algebra of orbifold models''},
Commun.\ Math.\ Phys.\ {\bf 123} (1989), 485--526.

\item{[DM]} {\smcap C.\ Dong and G.\ Mason}, 
{``Vertex operator algebras and Moonshine: a survey''}, 
In: Adv.\ Stud.\ Pure Math.\ {\bf 24} (Math.\
Soc.\ Japan, Tokyo, 1996), ed.\ by E.\ Bannai and A.\ Munemasa, pp.101--136.

\item{[DLM]} {\smcap C.\ Dong, H.\ Li, and G.\ Mason}, {``Modular-invariance of
trace functions in orbifold theory''}, preprint q-alg/9703016.

\item{[E]} {\smcap W.\ Eholzer}, {``On the classification of modular fusion algebras''},
Commun.\ Math.\ Phys.\ {\bf 172} (1995), 623--660.

\item{[ES]} {\smcap W.\ Eholzer and N.-P.\ Skoruppa},
{``Modular invariance and
uniqueness of conformal characters''}, Commun.\ Math.\ Phys.\ {\bf 174}
(1995), 117--136.

\item{[EZ]} {\smcap M.\ Eichler and D.\ Zagier}, {\it The Theory of Jacobi Forms},
(Birkh\"auser, Boston, 1985).

\item{[GG]} {\smcap M.\ Gaberdiel and T.\ Gannon}, {``The characters of rational conformal field theory
and vertex operator algebras''}, work in progress.

\item{[H]} {\smcap T.\ Hsu}, {``Identifying congruence subgroups of the modular 
group''}, Proc.\ Amer.\ Math.\ Soc.\ {\bf 124} (1996), 1351--1359.

\item{[H3]} {\smcap Y.-Z.\ Huang}, {``Vertex operator algebras and conformal
field theory''}, Int.\ J.\ Mod.\ Phys.\ {\bf A7} (1992), 2109--2151.

\item{[H2]} {\smcap J.\ Hurrelbrink}, {``On presentations of SL$_n(\Z_S)$''},
Commun.\ Alg.\ {\bf 11} (1983), 937--947.

\item{[J]} {\smcap G.\ A.\ Jones}, {``Congruence and non-congruence subgroups
of the modular group: a survey''}, In: {\it Proceedings of Groups---St.\
Andrews 1985} (Cambridge University, Cambridge, 1986), pp.223--234.

\item{[KP]} {\smcap V.\ G.\ Kac and D.\ Peterson}, {``Infinite--dimensional Lie algebras,
theta functions and modular forms''}, Adv.\ Math.\ {\bf 53} (1984), 125--264.

\item{[KSSB]} {\smcap T.\ H.\ Koornwinder, B.\ J.\ Schroers, J.\ K.\ Slinkerland, and
F.\ A.\ Bais}, {``Fourier Transform and the Verlinde Formula for the 
quantum-double of a finite group''}, math.QA/9904029.

\item{[L]} {\smcap S.\ Lang}, {\it Elliptic Functions}, 2nd edn.\ (Springer, 1987).

\item{[MD]} {\smcap D.\ L.\ McQuillan}, {``Some results on the linear fractional
groups''}, Illinois J.\ Math.\ {\bf 10} (1966), 24--38;

\item{} {\smcap J.\ B.\ Dennin, Jr.}, {``The automorphisms and conjugacy
classes of $LF(2,2^n)$''}, Illinois J.\ Math.\ {\bf 19} (1975), 542--552.

\item{[M]} {\smcap J.\ Mennicke}, 
{``On Ihara's modular group''}, Invent.\ math.\ {\bf 4} (1967), 202--228.

\item{[Mo]} {\smcap G.\ Moore}, {``Atkin-Lehner symmetry''}, Nucl.\ Phys.\ {\bf
B293} (1987), 139--188.

\item{[MS]} {\smcap G.\ Moore and N.\ Seiberg}, {``Classical and quantum conformal
field theory''}, Commun.\ Math.\ Phys.\ {\bf 123} (1989), 177--254.
  
\item{[O]} {\smcap A.\ Ogg}, {\it Modular forms and Dirichlet series},
(Benjamin, New York, 1969).

\item{[TNW]} {\smcap S.\ Tanaka}, {``Irreducible representations of the binary
modular congruence groups mod $p^\lambda$''}, J.\ Math.\ Kyoto Univ.\
{\bf 7} (1967), 123--132;

\item{} {\smcap A.\ Nobs and J.\ Wolfart}, {``Die irreduziblen Darstellungen der
Gruppen $SL\sb{2}(Z\sb{p})$, insbesondere $SL\sb{2}(Z\sb{2})$. I,II.''},
Comment.\ Math.\ Helvetici {\bf 51} (1976), 465--526.

\item{[Z]} {\smcap Y.\ Zhu}, {``Modular invariance of characters of vertex operator
algebras''}, J.\ Amer.\ Math.\ Soc.\ {\bf 9} (1996), 237--302.

\end